\documentclass[12pt]{article}
\usepackage{amsfonts}
\usepackage{graphics}
\usepackage{amsmath}

\newtheorem{Theorem}{Theorem}

\newtheorem{Corollary}[Theorem]{Corollary}

\newtheorem{Definition}[Theorem]{Definition}

\newtheorem{Lemma}[Theorem]{Lemma}

\newtheorem{Fundamental Theorem}{Fundamental Theorem}

\newenvironment{Proof}[1][Proof]{\textbf{#1.} }{\ \rule{0.5em}{0.5em}}

\def \fin {{\mathrm{fin}}}

\def \ad {\mathrm{ad}}
\def \tn {\mathrm{\otimes}}
\def \tr {\mathrm{tr}}
\def \Im {\mathrm{Im}}

\def \Jz  {\frac{J^z}{2z+1}}

\def \T {\mathcal{T}}

\def \Ct {\mathcal{C}}

\def \sgn {\mathrm{sgn}}
\def \e {\epsilon}
\def \a {\alpha}
\def \g {\mathfrak{g}}
\def \r {\rho}
\def \f {\phi}
\def \l {\lambda}
\def \s {\sigma}
\def \t {\tau}
\def \x {\xi}
\def \d {\delta}

\def \h {[[h]]}
\def \ei {{e^{i\theta}}}

\def \mm {{\underline{m}}}
\def \D {\Delta}
\def \X {\chi}
\def \W {\Omega}
\def \A {\mathcal{A}}
\def \C {\mathbb{C}}
\def \Q {\mathrm{Q}}
\def \R {\mathbb{R}}

\def \Z {\mathcal{Z}}
\def \L {\mathcal{L}}
\def \Int {\mathbb{Z}}
\def \N {\mathbb{N}}
\def \L {\mathcal{L}}
\def \B {\mathcal{B}}
\def \S {\mathcal{S}}
\def \sl {\mathfrak{sl}(2,\C)}
\def \cZ {\mathcal{Z}}

\begin{document}

\title{On the Analytic Properties of the $z$-Coloured Jones Polynomial}
\author{Jo\~{a}o  Faria Martins\footnote{email address: pmxjm@maths.nottingham.ac.uk}}

\date{\today}
\maketitle

\begin{abstract} We analyse the possibility of defining $\C$-valued Knot invariants associated with infinite dimensional unitary representations of $SL(2,\R)$ and the Lorentz Group taking as starting point the Kontsevich Integral and the notion of infinitesimal character. This yields a family of knot invariants whose target space is the set of formal power series in $\C$, which contained in the Melvin-Morton expansion of the coloured Jones polynomial. We verify that for some knots the series have zero radius of convergence and analyse the construction of functions of which this series are asymptotic expansions by means of Borel re-summation. Explicit calculations are done in the case of torus knots which realise an analytic extension of the values of the coloured Jones polynomial to complex spins. We present a partial answer in the general case.
\end{abstract}

\section*{Introduction}

Since the advent of quantum groups, and in particular of the notion of a 
quantised universal enveloping algebra of a semisimple Lie algebra in the end of the eighties, their theory has been applied to the construction of link
invariants. The main idea behind all approaches comes from  the observation
that, in the current terminology, they are ribbon Hopf algebras \cite{RT}, 
which implies that their category of finite dimensional representations is a ribbon category, with trivial associativity constraints. In the pioneering work of Freyd and Yetter, cf. \cite{FY}, it was observed that the (ribbon) tangles form a ribbon category which is universal in the class of all strict ribbon categories. This framework gives us a knot invariant for any ribbon Hopf algebra and any finite dimensional representation of it, as observed in the construction of Reshetikin and Turaev's functor defined in \cite{RT}.

A limitation of the constructions above is that they are not directly 
applicable to the case of infinite dimensional representations of ribbon Hopf algebras. This is because they involve taking traces or the use of coevaluation maps, which are difficult to define in the infinite dimensional context. 
However, one is forced to consider knot invariants associated with infinite dimensional representations when generalising to invariants associated with unitary representations of non-compact groups. This kind of representation appears in the context of $(2+1)$-quantum gravity and Chern-Simons theory with non-compact groups. See for example \cite{W},\cite{BC},\cite{BNR}, \cite{GI}, \cite{NR} or \cite{G}. It would thus be important  to define $\C$-valued knot invariants associated with representations of this kind. The aim of this paper, which continues \cite{FM}, is describe a possible path for doing this. We shall be mostly interested in the $SL(2,\R)$ and $SL(2,\C)$ cases.

The $h$-adic quantised universal enveloping algebras of semisimple Lie algebras are usually easier to deal with in the context of infinite dimensional representations. In this article we shall restrict to them. In this case there are various different variants of the construction of quantum Knot invariants. Some of them can be used in the infinite dimensional case, for example,  Laurence's Universal $U_h(\g)$ knot invariant or the Kontsevich Universal knot invariant. Roughly speaking, given a (complex semisimple) Lie algebra and an $\ad$-invariant non-degenerate symmetric bilinear form on it, they will yield a knot invariant which take values in the algebra of formal power series over the centre of $U(\g)$, the universal enveloping algebra of $\g$. It is called the universal $U(\g)$-knot invariant.   

The main idea behind the construction of non-compact group knot invariants is the following. Suppose we have a representation $\r$ of the Lie algebra $\g$ in some vector space $V$, which we do not assume to be finite dimensional. We can always lift it to a representation, also denoted by $\r$, of the enveloping algebra of $\g$. In some cases it can happen that any element of the centre of $U(\g)$ acts in $V$ as a multiple of the identity. Such representations thus define  an algebra morphism from the centre of $U(\g)$ to $\C$, a  central character of $U(\g)$. They are usually called representations  which admit a central character. This type of $\g$-module arises naturally in Lie algebra theory. Some examples would be the cyclic highest weight representations of a semisimple Lie algebra. Notice they are infinite dimensional if the weight is not integral. It is a well established fact that the central characters of them exhaust all central characters of $U(\g)$ if $\g$ is complex semisimple.
Another context where representations which admit  a central character appears is the context of unitary irreducible representations $R$ of real Lie Groups $G$ in complex Hilbert spaces $V$. More precisely, it is possible to prove that the induced representation $R^{\infty}$ of $U(\g\tn_\R \C)$  in the space of smooth vectors of $V$ under the action of $R$ admits a  central character. Here $\g$ denotes the Lie algebra of $G$. It is called the infinitesimal character of $R$.

Any central character of $U(\g)$ can be used to evaluate the universal $U(\g)$ knot invariant. This will then yield a knot invariant with values in the space of formal power series over $\C$.  Obviously, one price we have to pay when we  consider infinite dimensional representation of $\g$  will then be, in general, the need to stick to representations of $\g$ that admit a central character  and links with one component (knots). In this article we propose to consider this kind of knot invariant in the context of irreducible unitary representations of $SL(2,\R)$ and $SL(2,\C)$. Notice that as they are non-compact Lie groups they admit infinite dimensional irreducible unitary representations. 

As mentioned before, in  the semisimple Lie algebras context  the central characters of the highest weight representations exhaust all central characters of $U(\g)$. Moreover the value of these central characters in a central element of $U(\g)$ depends polynomially on the weight and it is  determined by its values on the weights that define finite dimensional representations. In particular the  knot invariants obtained by admitting infinite dimensional representations are not more powerful than the already known ones. In the power series level they are in an obvious sense  analytic continuation of the usual quantum groups knot invariants. The fact finite dimensional representations suffice was also pointed out in \cite{BN}. For example in the $SL(2,\R)$ and $SL(2,\C)$ context, these non compact knot invariants express out of a analytic continuation of  the coloured Jones polynomial to complex spins, termwise in  the power series. In the former case  this extension is immediate from the Melvin-Morton expansion of it. We called this extension the $z$-coloured Jones polynomial. We need however an infinite set of finite representations to determine these knot invariants. It is unclear what happens in the non-semisimple Lie algebras context.

As we referred before, the definition of $\C$-valued, that is numerical, knot invariants would be the most important for applications. So we want to say something about about the kind power series that we obtain. This will be one of the main subjects of this paper. A main result will be that for a large class of interesting unitary infinite dimensional representations of $SL(2,\R)$ the associated series has a zero radius of convergence, at least in the case of torus knots. The same is true in the $SL(2,\C)$ case. Notice that this does not happen in the case of finite dimensional representations, the case in which there is no problem in defining numerical knot invariants. However, in the context of torus knots, they define Borel summable series. This means there is a natural way to find analytic functions of which these power series are asymptotic expansions. Also, that the uncertainty in process of re-summation is reduced to a numerable, in this case finite, set of functions, differing by rapidly decreasing terms. It would be interesting to analyse what this uncertainty means. This re-summation realises an analytic extension of  the coloured Jones polynomial of torus knots to complex spins, in the context of numerical knot invariants rather than only termwise in the  power series,  that is, of the actual values of the coloured Jones polynomial. In the general case of an arbitrary knot and a unitary representation of $SL(2,\R)$ or $SL(2,\C)$ it is possible to prove that the series obtained are of Gevrey type $1$. This is a necessary condition for Borel summability and permits us to define a weaker process of re-summation up to exponentially decreasing functions. It is an open problem whether the process of Borel re-summation of the $z$-coloured Jones polynomial works for any knot.

\section{Preliminaries} 

\begin{figure}
\includegraphics{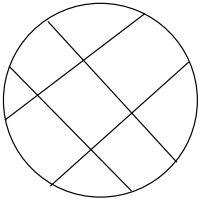}
\caption{\label{chord} A chord diagram with $4$ chords}
\end{figure}

We recall the definition  of the algebra of  chord diagrams, which is the
target space for the Kontsevich Universal Knot Invariant. For more details see
for example \cite {BN} or \cite{K}. Both references contain almost all the material considered in this section. 
A chord diagram is a finite set  $w=\{c_1,...,c_n\}$ of non intersecting, cardinality 2 subsets of the oriented circle, modulo orientation preserving homeomorphisms. The subsets $c_k$ are called chords and are supposed to be pairwise disjoint. We usually specify a chord diagram by drawing it as in figure \ref{chord}. In all the pictures we assume the circle is  oriented counterclockwise.

\begin{figure}
\includegraphics{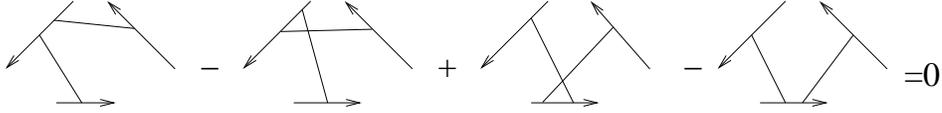}
\caption{\label{4T} $4$ Term relations}
\end{figure}

For each $n \geq 2$, let $V_n$ be the free $\C$ vector space on the set of all
chord diagrams with $n$ chords. That is the set of formal finite linear
combinations $w=\sum_i \l_i w_i$, where $\l_i \in \C$ and $w_i \in V_n$ for
any $i$. Consider the sub vector space $4T_n$ of $V_n$ which is the subspace
generated by all linear combinations of chord diagrams of the form displayed
in figure \ref{4T}. The three intervals considered in the circle can appear at
an arbitrary order in $S^1$. Define for each $n \in \{0,1,2,..\}$, the vector
space $\A_n=V_n/4T_n$. We consider $\A_0=V_0$ and $\A_1 =V_1$.

Let $\g$ be a Lie algebra over the field $\C$. An infinitesimal R-matrix, cf \cite{K}, in 
$\g$ is a symmetric tensor $t \in \g \tn \g$ such that $[\D(X),t]=0,
\forall X \in \g$. The commutator is taken in  $U(\g) \tn U(\g)$, where
$U(\g)$ denotes the universal enveloping algebra of $\g$. The map
$\D:U(\g) \to U(\g) \tn U(\g)$ is the standard coproduct in
$U(\g)$. It verifies $\D(X)=X\tn 1 +1\tn X$ if $ X \in \g$.

Suppose we are given an infinitesimal R-matrix $t$. Write $t=\sum_i a_i \tn
b_i$ with $a_i,b_i \in \g$. We  have $\sum_j [\D(a_j),t]\tn b_j=0$, thus
\begin{equation}
\sum_{i,j} a_j a_i \tn b_i \tn b_j- a_i a_j \tn b_i \tn b_j+ a_i \tn a_j b_i \tn b_j- a_i \tn b_i a_j \tn b_j=0,
\end{equation}
which resembles the $4T$ relations just considered. Given a chord
diagram $w$ and an infinitesimal R-matrix $t=\sum_i a_i \tn b_i$ it is thus natural to construct an element $\f_t(w)$ of $U(\g)$ in the following fashion, cf \cite{K}: Start in an arbitrary point of the circle and go around it in the direction of its orientation. Order the chords of $w$ by the order with which you pass them as
in figure \ref{chords}. Each chord has thus an initial and an end
point. Then go around the circle again and write (from the right to the left)
$a_{i_k}$ or $b_{i_k}$ depending on whether you got to the initial or final
point of the $k^{th}$ chord. Finally, sum over all the $i_k$'s. For example for
the chord diagram of picture  \ref{chords} the element $\f_t(w)$ is:
\begin{equation}
\sum_{i_1,i_2,i_3} b_{i_2} b_{i_3} b_{i_1} a_{i_3} a_{i_2} a_{i_1}.
\end{equation}
Let $w$ be a chord diagram. A  good feature about the element $\f_t$ is:
         
\begin{figure}
\includegraphics{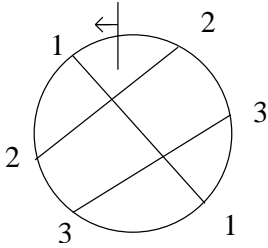}
\caption{\label{chords} Enumerating the Chords of a Chord Diagram}
\end{figure}
  
\begin{Theorem} \label{central} Let $\g$ be a complex Lie algebra and $t \in \g \tn \g$ be an infinitesimal R-matrix of $\g$. Let also $w$ be a chord diagram with $n$-chords. The element $\f_t(w)$ is well defined, that is it does not depend on the chosen point of the chord diagram $w$, and it is a central element of $U(\g)$. Moreover the assignment $w \mapsto \f_t(w)$ descends to a linear map $\f_t:\A_n \to \Ct(U(\g))$. Here $\Ct(U(\g))$ denotes the centre of $U(\g)$.
\end{Theorem}
See for example \cite{K} or \cite{CV}.

In general infinitesimal R-matrices in a Lie algebra $\g$ are constructed out of a non-degenerate, $\ad$-invariant symmetric bilinear  form $<,>$ in $\g$. Invariance here means $<[X,Y],Z]>+<Y,[X,Z]>=0,\forall X,Y,Z \in \g$. Concretely let $\{X_i\}$ and $\{X^i\}$ be  a  basis and a dual basis of $\g$ with respect to $<,>$. It is easy to show that the tensor $\sum_i X_i \tn X^i$ is an infinitesimal R-matrix of $\g$. See also \cite{CV}.

The details of the definition of the Kontsevich Integral can be easily found in the literature. The classical reference is \cite{BN}. For the framed Kontsevich Integral which we are going to use  see \cite{K} or \cite{LM}. We give a small review in the appendix. Recall that the Kontsevich Integral $\Z$ is a framed Knot invariant of the form 
\begin{equation}
\Z:K \mapsto  \sum_{n=0} ^{+\infty} \Z_n(K) h^n,
\end{equation}
where $\Z_n(K) \in \A_n, n=0,1...$ and $h$ is a formal variable. We shall take the normalisation of $\Z$ for which the value on the unknot is the Wheels element, cf \cite{LNT}. Let $\g$ be a Lie algebra and $t$ an infinitesimal $R$-matrix of $\g$. We can consider the composition $\f_t \circ \Z$. This will yield a knot invariant with values in the algebra of formal power series in the centre of $U(\g)$. Therefore:
\begin{Theorem} \label{universalg}
Let $\g$ be Lie algebra and $t$ an infinitesimal R-matrix in $\g$. There exists a framed knot invariant $(\f_t \circ \Z)$. It has the form: 
\begin{equation}
(\f_t\circ \Z): K \mapsto \sum_{n=0}^{+\infty} (\f_t \circ \Z_n)(K)h^n, 
\end{equation}
where $(\f_t \circ \Z)(K) \in \Ct(U(\g)),n=0,1..$.
\end{Theorem}

If $\g$ is complex semisimple and $t$ is the infinitesimal R-matrix coming from its Cartan-Killing form in $\g$, then $(\f_t \circ \Z)(K)$ defines an analytic function from $\C$ to a completion of $U(\g)$, where $U(\g)$ is given the topology of convergence in its finite dimensional representations, cf. \cite{PS}. In other words, for any finite dimensional representation $\r$ of $\g$ in $V$, the power series   $\r(\f_t \circ\Z)(K)$ converges to a linear operator $V \to V$, in fact to a multiple of the identity,  if $V$ is irreducible. The usual quantum group knot invariants are obtained by taking the trace of these operators. We will go back to these issues later. This article aims mostly to consider the case in which we admit infinite dimensional representations.

\section{Non Compact Group Knot Invariants}\label{NONCOMPACT}

Let $G$ be a Lie group, always assumed to be real, and $\g$ its Lie algebra. Let also $\g_\C=\g\tn_\R \C$ denote the complexification of $\g$. Consider a unitary representation $R$ of $G$ in the complex Hilbert space $V$. Notice $V$ is not assumed to be finite dimensional. We recall that unitarity means that the operator $R(g):V \to V$ is a unitary for any $g \in G$. Also we suppose a continuity condition, namely that for any $v \in V$ the map $g \in G \mapsto R(g)v \in V$ is continuous.

 The main reference for what follows is \cite{Kir}. Let $V_\infty$ denote the vector space of smooth vectors of $R$. That is 
\begin{equation}
V_\infty = \{v \in V: g \in G \mapsto R(g)(v) \in C^{\infty}(G,V) \}.
\end{equation} 
It is well known $V_{\infty}$ is dense in $V$. Differentiating $R$ at the identity of $G$  defines a map $R^{\infty}:\g \tn V_{\infty} \to V$. It is possible to show that $V_{\infty}$ is invariant under $\g$ and that $R^{\infty}$ is a honest representation of $\g$. It extends therefore to a representation, which we also call $R^{\infty}$, of $U(\g_\C)$ in $V_\infty$. Recall that $V$ is a complex vector space.

Suppose $R$ is a irreducible representation. In our context this means that $V$ has no closed invariant subspaces under the action of $G$. 
The following (non trivial) result can be found in \cite{Kir}:

\begin{Theorem}\label{infinitesimal}
If $V$ is irreducible then any element of $\Ct(U(\g_\C))$ acts on $V_{\infty}$ under $R^\infty$ as a multiple of the identity operator. 
\end{Theorem}

Recall  $\Ct(U(\g_\C))$ denotes the centre of the universal enveloping algebra of $\g_\C$. Obviously if $V$ is a finite dimensional complex vector space then the theorem just stated is a consequence of Schur's lemma and the unitarity condition is not needed. Also it is in general possible to show directly that the above property is true for a large class of infinite dimensional representations, not necessarily unitary. Amongst them are the representations of $SL(2,\R)$ and $SL(2,\C)$ which we are going to consider. However the last theorem tells us that our construction is general. 
  
Another way to state theorem \ref{infinitesimal} is to say that $R$ has an infinitesimal character. In other words there exists a (unique) central character  $\X_R^\infty$ of $U(\g_\C)$, that is a morphism of complex algebras $\Ct(U(\g_\C)) \to \C$, with the property: 
\begin{equation}
R^{\infty}(a)(v)=\X_R^{\infty}(a)v,\forall v \in V^{\infty},\forall a \in \Ct(U(\g_\C)).
\end{equation} 

Therefore we have the following obvious conclusion:
\begin{Theorem}\label{defofI}
Let $G$ be a real Lie group and $R$ an irreducible unitary representation of $G$ in some complex Hilbert space $V$. Let also $t$ be an infinitesimal R-matrix in the complexification $\g_\C$ of the Lie algebra  $\g$ of $G$. There exists a knot invariant with values in $\C[[h]]$:
\begin{equation}
I(G,t,R)=\X^{\infty}_R \circ \f_t \circ \Z.
\end{equation}
It has the form
\begin{equation}
K \mapsto \sum_{n=0}^{+\infty} \X_R^{\infty}(\f_t \circ \Z_n(K))h^n.
\end{equation}
\end{Theorem}
As an example, consider  $G=SU(2)$, thus $\g=\mathfrak{su}(2)$ and $\g_\C=\mathfrak{sl}(2,\C)$. Take  $t$ to be the infinitesimal R-matrix coming from minus the Cartan-Killing form in $\mathfrak{sl}(2,\C)$, that is $<X,Y>=-\tr( \ad(X) \circ \ad(Y))$. The tensor $t$ has the form:
\begin{equation}
t=-\frac{1}{4}(\s_X\tn \s_X+\s_Y\tn \s_Y+\s_Z\tn \s_Z),
\end{equation}
where 
\begin{equation}
\s_X=\frac {1}{2}\begin{pmatrix} i & 0 \\ 0 & -i \end{pmatrix},\s_Y=\frac
{1}{2}\begin{pmatrix} 0 & i \\ i & 0 \end{pmatrix},\s_Z=\frac
{1}{2}\begin{pmatrix} 0 & -1 \\ 1 & 0 \end{pmatrix}.
\end{equation}
The following is well known. 
\begin{Theorem}\label{finite}
Let $J^{\a}$ denote the framed coloured Jones polynomial associated with the representation $R^\a$ of $SU(2)$ with  spin $\a\in \{0,1/2,1,3/2...\}$. We take the normalisation of the coloured Jones polynomial that evaluates at the unknot to the quantum dimension of the  $U_h(\mathfrak{sl}(2,\C))$ spin $\a$ representation. Given any framed knot $K$ we have:
\begin{equation}
\frac{J^{\a}}{2\a+1}(K)=I(SU(2),t,R^\a)(K),\a=0,1/2,1,3/2....
\end{equation}
as formal power series.
\end{Theorem}
This  a non trivial result. A  path for proving it relies upon the  framework of quasi Hopf algebras and the notion of gauge transformations on them. This is  described by Drinfeld in \cite{D1} and \cite{D2}. The context in which we need to apply it is the  quantised universal enveloping algebras one. In this case some rigidity results ensure the above theorem is true. All this is  described in detail in the same references. For a complete discussion, see  \cite{LM} or \cite{K}.

\subsection{Some Examples in the $SL(2,\R)$ case}\label{repsofSL}

Let us now describe some infinite dimensional examples. Consider the Lie group $G=SL(2,\R)$. It is a non-compact semisimple group. As before we have $\g_\C=\mathfrak{sl}(2,\C)$. Take again $t$ to be the infinitesimal R-matrix coming from minus the Cartan-Killing form in $\g_\C$. We can also write it as 

\begin{equation}
t=-\frac{1}{4}\left (E \tn F +F \tn E +\frac{H \tn H}{2} \right ),
\end{equation}
where:
\begin{equation}
H=\begin{pmatrix}1 & 0 \\ 0 & -1 \end{pmatrix},  E=\begin{pmatrix}0 & 1 \\ 0
   & 0 \end{pmatrix}, F=\begin{pmatrix}0 & 0 \\ 1 & 0 \end{pmatrix}.
\end{equation}
We start by defining the representations of $SL(2,\R)$ in the principal series. They depend on a imaginary parameter $s \in i\R$ and an $\e \in \{0,1\}$, the parity of the representation. In general $V^{s,\e}$ is the space $L^2(\R)$ of complex-valued square integrable functions in $\R$. The action of $SL(2,\R)$ has the form:
\begin{equation}
R^{s,\e}\left (\begin{pmatrix}a & b \\ c & d \end{pmatrix}
\right)(f)(x)=\sgn^\e(bx+d)|b x +d|^ {s-1}f\left(\frac{a x+c}{b x+d} \right ).
\end{equation}
See \cite{GI} for a alternative description of these representations, as well as the definition of their associated  spin network theory. It applies to the construction of spin foam models for $(2+1)$-Quantum Gravity.

The positive discrete series depends on a parameter $m \in \Int^-$. The representations of this type are denoted by $R^{m,+}$. The representation space $V$ for the representation $\R^{m,+}$ is the space of holomorphic functions $f$ in the upper half plane such that:
\begin{equation}
\frac{i}{2\Gamma(-m)}\int_{\Im(z)>0}|f(z)|^2(\Im(z))^{-m-1}d z d \bar{z} <+\infty.
\end{equation}
The inner product in $V$ as an expression similar to the formula above. The group $SL(2,\R)$ acts in the fashion:
\begin{equation}
R^{s,\e}\left (\begin{pmatrix}a & b \\ c & d \end{pmatrix} \right)(f)(x)=(b x +d)^{s-1}f\left(\frac{a x+c}{b x+d} \right ).
\end{equation}

The representations in the negative series are denoted by $R^{m,-}$. They depend on a parameter $m \in \Int^-$. The representation space for them is the space of holomorphic functions $f$ in the lower half plane such that:
\begin{equation}
\frac{1}{2\Gamma(-m)}\int_{\Im(z)>0}|f(z)|^2|\Im(z)|^{-m-1}d z d \bar{z} <+\infty.
\end{equation}
The action of $SL(2,\R)$ in $V$ is similar to the case of the positive discrete series.

The Lie group $SL(2,\R)$ still has one more series of unitary representation, namely the complementary series of representations. Details can be found in \cite{L} or \cite{GGV}. 
Excluding the representation $R^{0,1}$, all the representations considered are unitary and irreducible. Therefore there is attached to them a knot invariant with values in $\C[[h]]$. In the next section we relate these knot invariants to the coloured Jones polynomial.

\subsection{The z-coloured Jones polynomial}\label{zcoloured}

Let $G$ be a real Lie group, $\g$ its Lie algebra and $\g_\C$ the complexification of $\g$. We suppose $\g_\C$ is equipped with a infinitesimal R-matrix $t\in \g_\C \tn \g_\C$. Let also $R$ be a irreducible unitary representation of $G$ in complex Hilbert space $V$. A closer look at the definition of the knot invariant $I(G,t,R)$, cf theorem \ref{defofI},  tells us that the only information which we took from the representation  $R$ was its infinitesimal character $\X_{R}^{\infty}$. It is a morphism of algebras from the centre $\Ct(U(\g_\C))$ of $U(\g_\C)$ to $\C$, that is a central character of $U(\g_\C)$. In the case $\g_\C$ is semisimple, we know the form of all such morphisms. Let us say what the situation is in the case $\g_C=\mathfrak{sl}(2,\C)$.  We refer to \cite{V} for further details. All this generalises for any semisimple Lie algebra. In particular any similar construction of knot invariants out of infinite dimensional representations of semisimple Lie groups will have the same kind of properties.  

Consider any Cartan decomposition of $\mathfrak{sl}(2,\C)$ and choose a Borel subalgebra relative to it. Given a complex number $a$ there exists a unique cyclic, highest weight representation $\r^{\frac{a}{2}}$ of maximal weight $a$. This representation is finite dimensional if, and only if, $a =0,1,2,...$. In this case it integrates to the representation of $SL(2,\C)$ of spin $\frac{a}{2}$. If $a$ is an arbitrary complex number, any element of the centre of $U(\mathfrak{sl}(2,\C))$ acts in the representation space of $\r^{\frac{a}{2}}$ as a multiple of the identity. Denote by $\l_{\r^{\frac{a}{2}}}$ the central character of $\r^{\frac{a}{2}}$. That is if $x$ is a central element of $U(\mathfrak{sl}(2,\C))$ then $\l_{\r^{\frac{a}{2}}}(x)$ is the unique complex number such that 
$\l_{\r^{\frac{a}{2}}}(x)v=\r^{\frac{a}{2}}(x)(v)$ for all $v$ in the representation space of $\r^{\frac{a}{2}}$. We will then have $\l_{\r^{\frac{a}{2}}}=\X_{R^{\frac{a}{2}}}^{\infty}, a\in\{0,1,2,..\}$. Recall $R^\a$ denotes the representation of $SU(2)$ of spin $\a=0,1/2,1,3/2...$. 

All the pieces of the following theorem can be found in \cite{V}.

\begin{Theorem} \label{character}
Let $f:\Ct(U(\mathfrak{sl}(2,\C))) \to \C$ be a morphism of complex algebras. In other words, a central character of $U(\mathfrak{sl}(2,\C))$. We have:
\begin{enumerate}

 \item There exists an $a \in \C$ such that $f=\l_{\r^{\frac{a}{2}}}$.
 \item $\l_{\r^{\frac{a}{2}}}=\l_{\r^{\frac{b}{2}}}$ if and only if $(a+1)^2=(b+1)^2$.
 \item Given  $x \in \Ct(U(\mathfrak{sl}(2,\C)))$ the map $a\in \C  \mapsto \l_{\r^{\frac{a}{2}}}(x)$ is a polynomial in $a$ of degree smaller or equal to the degree of $x$ in $U(\mathfrak{sl}(2,\C))$. In fact it is a polynomial in  $(a+1)^2$  
\end{enumerate} 
\end{Theorem}
Given $z\in \C$, it is therefore natural to define a $z$-coloured Jones polynomial in the form: (cf theorem \ref{finite})
\begin{equation}
\frac{J^z}{2z+1}(K)=\l_{\r^{z}} \circ (\f_t \circ \Z)(K),
\end{equation}
here $K$ denotes a framed knot. Due to part $3.$ of Theorem \ref{character} we have:
\begin{equation}\label{melvinmorton}
\frac{J^{z}}{2z+1}(K)=\sum_{n=0}^{+\infty}J_n^z(K)h^n=\sum_{n=0}^{+\infty}\left( \sum_{k=0}^{2n} J_{n,k}(K) z^k \right)h^n
\end{equation}
Notice the degree of $(\f_t \circ \Z_n)(K)$ in $U(\mathfrak{sl}(2,\C)$ is not bigger than $2n$, for any framed knot $K$. Obviously equation (\ref{melvinmorton}) is exactly the Melvin-Morton expansion of the coloured Jones polynomial for a given knot $K$, cf \cite{MM} or \cite{GN}.  

Some  properties of the $z$-coloured Jones polynomial coming from the corresponding properties of the Melvin-Morton expansion are the following:

\begin{enumerate}
\item If $K$ is a framed knot then $J^z_n(K)$ is a polynomial in $(2z+1)^2$ of degree smaller or equal to $n$. 
\item If $(2z+1)^2=(2w+1)^2$ then $\frac{J^z}{2z+1}=\frac{J^w}{2w+1}$.

\item $\frac{J^z}{2z+1}$ is the usual (rescaled) Jones polynomial if $2z+1=1,2,...$.
\item  If $2z+1=1,-1,2,-2...$ then $\frac{J^z}{2z+1}(K)$ defines a power series in $h$ convergent in $\C$.
\end{enumerate}
Properties $1$, $2$ and $3$ are easy consequences of our discussion. The
fourth is a consequence of the fact that, in our normalisation, $\frac{J^z}{2z+1}(K)$ is a Laurent
polynomial in $e^{h/4}$ if $2z+1=1,-1,2,-2...$, see \cite{MM}. This generalises to any semisimple Lie algebra. 

\subsubsection{Some examples}\label{someexamples}
Let $q=\exp(h)$. For an $n\in \N_0$ and a $z \in \C$, consider the term 
\begin{equation}
D(n,z)=\prod_{k=1}^{n}\left[ \left(q^{\frac{2z+1}{2}}-q^{-\frac{2z+1}{2}}\right ) -\left(q^{\frac{k}{2}}-q^{-\frac{k}{2}}\right )\right ].
\end{equation} 
It is a power series in $h$ such that the first $2n$ terms are zero. Therefore
if $f(n) ,n \in \N_0$  are power series in $h$, for example  Laurent
polynomials in $q$ and $q^{-1}$, then $\sum_{n \in \N_0} f(n)D(n,z)$ is an
infinite series of power series which is termwise convergent, since it is of
terminating type, termwise. Suppose
$A(z)=\sum_{n\in \N_0} A_n(z)h^n$ and $B(z)=\sum_{n \in \N_0} B_n(z) h^n$ are power series whose coefficients depend polynomially in $z$, for example power series such as $q^{\frac{2z+1}{2}}$ or $q^{-\frac{2z+1}{2}}$. Then also the coefficients of their product
depend polynomially in $z$, thus in particular  the coefficients of $D(n,z)$, for any $n \in \N_0$. The same is true for the coefficients of any power series of the type $\sum_{n \in \N_0} f(n)D(n,z)$,  where  $f(n) \in \C \h$ is such that its terms depend polynomially in  $z$.

Let $3_1$ and $4_1$ denote the zero framed trefoil and figure of eight knots. We have, see \cite{H}:
\begin{equation}
\frac{J^{z}}{2z+1}(3_1)=\frac{1}{2z+1}\frac{q^{\frac{2z+1}{2}}-q^{-\frac{2z+1}{2}}}{q^{\frac{1}{2}}-q^{-\frac{1}{2}}}\sum_{n=0}^{+\infty} (-1)^n q^{-n(n+3)/2}D(n,z),
\end{equation}
\begin{equation}
\frac{J^{z}}{2z+1}(4_1)=\frac{1}{2z+1}\frac{q^{\frac{2z+1}{2}}-q^{-\frac{2z+1}{2}}}{q^{\frac{1}{2}}-q^{-\frac{1}{2}}}\sum_{n=0}^{+\infty}f_n D(n,z).
\end{equation} 
In general for any framed knot $K$ there exist Laurent polynomials $f_n(K)(h), n
\in \N_0$ in  $q$ and $q^{-1}$ such that 
\begin{equation}\label{Habiro}
\frac{J^{z}}{2z+1}(K)=\frac{q^{F(K)\frac{z(z+1)}{2}}} {2z+1}\frac{q^{\frac{2z+1}{2}}-q^{-\frac{2z+1}{2}}}{q^{\frac{1}{2}}-q^{-\frac{1}{2}}}\sum_{n=0}^{+\infty}f_n(K),
D(n,z) 
\end{equation}
where $F(K)$ is the framing coefficient of $K$. Actually, Habiro proved these formulae only in the case of finite
 dimensional representations, that for $z \in \frac{1}{2} \N_0$. However the coefficients of the power series
 above  depend polynomially in $z$, which implies the formulae are true also
 for infinite dimensional representations. This is the old principle that if
 two polynomials coincide in an infinite set then they are the same. We shall
 use this method  of proof quite frequently. Equation \ref{Habiro} also
 proves that $\frac{J^{\a}}{2\a+1}(K)$ always defines a Laurent polynomial in $q^{1/4}$ if $\a =0,\frac{1}{2},1,...$.

\subsubsection{Back to $SL(2,\R)$}
 Recall the framed knot invariants $I(G,t,R)$ defined in theorem \ref{defofI}.
Given that $\sl$ is simple, it is easy to prove that any infinitesimal R-matrix in $\mathfrak{sl}(2,\C)$ is a multiple of the one coming from minus the Cartan-Killing form considered previously. The following result is a straightforward consequence of the discussion above: 

\begin{Theorem}
Let $G$ be a real form of $SL(2,\C)$ and $R$ be an irreducible unitary
representation of $G$ in the Hilbert space $V$. Let also  $t$ be an
infinitesimal R-matrix in $\mathfrak{sl}(2,\C)$. After rescaling $t$ (possibly), there exists a $z \in \C$
such that $I(G,t,R)=\frac{J^z}{2z+1}$.
\end{Theorem}  
 We can prove similar results for any semisimple Lie group.

It is not difficult to find the exact relation between the invariants associated with the unitary representations of $SL(2,\R)$  and the $z$-coloured Jones polynomial which needs to exist in the light of the theorem above. Let as usual $t$ be the infinitesimal R-matrix in $\mathfrak{sl}(2,\C)$ coming from minus its Cartan-Killing form. We refer to \cite{L} (chapter VI), for an infinitesimal description of the unitary representations of $SL(2,\R)$, in terms of a Hilbert basis of the representations spaces made out of analytic vectors. In general the highest weight representations $\r^{a}$ of $\mathfrak{sl}(2,\C)$ cannot be integrated to representations of $SL(2,\R)$, however it is possible to relate their central characters with the  infinitesimal characters of the unitary representations of $SL(2,\R)$, given that the second ones also  depend polynomially in the parameters defining them, as a glance in the infinitesimal expression of the representations tells us. In fact by an argument very similar to the proof of thorem $12$ of \cite{FM} we can prove: 

\begin{Theorem}\label{Continuation}
Let $t$ be the infinitesimal R-matrix in $\mathfrak{sl}(2,\C)$ coming from minus its Cartan-Killing form. We have, for the principal series of representations:
\begin{equation}
I(SL(2,\R),t,R^{s,\e})=\frac{J^{\frac{s-1}{2}}}{s}, s \in i \R \textrm{ and } \e\in \{0,1\}
\end{equation}
and
\begin{equation}
I(SL(2,\R),t,R^{m,\pm})=\frac{J^{\frac{m-1}{2}}}{m}, m \in \Int^-,
\end{equation}
for the discrete series.
\end{Theorem} 
A full proof of this fact will appear elsewhere.  
Observe that as a consequence the usual coloured Jones polynomial can be obtained out of the unitary infinite dimensional representations of $SL(2,\R)$ in the discrete series. 

Another consequence is the fact observed in the introduction that the quantum knot invariants associated with the infinite dimensional representations $SL(2,\R)$ are in a sense analytic continuations of the ones associated with finite dimensional representations of the complexification $\mathfrak{sl}(2,\C)$ of $\mathfrak{sl}(2,\R)$. Also the fact that the non-compact knot invariants are not stronger than the ones associated with finite dimensional representations. The meaning of this is obvious from the Melvin-Morton expansion (\ref{melvinmorton}). Notice however that this property is valid only termwise in the power series expansions, and we will see later that in some cases the analytically continued power series may have a zero radius of convergence. This does not happen in the case of the coloured Jones polynomial. We will come back to this at the end of this article. This property concerning analytic continuations is valid if $G$ is any semisimple Lie group. I do not know what the answer is in the case of non semisimple Lie algebras. 

\subsection{Lorentz Group Case}

We now pass to the description of the corresponding $SL(2,\C)$ invariants. We always look at $SL(2,\C)$ as a real Lie group. It is the universal covering of the Lorentz group. More details of the following description can be found in \cite{FM}. The Lorentz Lie algebra $L$ is defined as being the realification of the complex Lie algebra $\mathfrak{sl}(2,\C)$. A real basis of $L$ is given by 
\begin{equation}
\{\s_X,B_X=-i\s_X,\s_Y,B_Y=-i\s_Y,\s_Z,B_Z=-i\s_Z\},
\end{equation}
where, as before
\begin{equation}
\s_X=\frac {1}{2}\begin{pmatrix} i & 0 \\ 0 & -i \end{pmatrix},\s_Y=\frac
{1}{2}\begin{pmatrix} 0 & i \\ i & 0 \end{pmatrix},\s_Z=\frac
{1}{2}\begin{pmatrix} 0 & -1 \\ 1 & 0 \end{pmatrix}.
\end{equation}
We denote by $L_\C=L \tn _{\R} \C$ the complexification of $L$.  
There exists a unique isomorphism of complex Lie algebras $\t:\mathfrak{sl}(2,\C) \oplus \mathfrak{sl}(2,\C) \to L_\C$ such that 
\begin{equation}
\t(\s_X\oplus 0)=\frac{\s_X+iB_X}{2},\quad \t(0 \oplus \s_X)=\frac{\s_X-iB_X}{2},\end{equation}
\begin{equation}
\t(\s_Y\oplus 0)=\frac{\s_Y+iB_Y}{2},\quad \t(0 \oplus \s_Y)=\frac{\s_Y-iB_Y}{2},\end{equation}
\begin{equation}
\t(\s_X\oplus 0)=\frac{\s_Z+iB_Z}{2},\quad \t(0 \oplus \s_Z)=\frac{\s_Z-iB_Z}{2}.
\end{equation}
In general we denote  $x^l=\t(x \oplus 0)$ and $x^r=\t(0 \oplus x)$. The notation $a^l$ and $a^r$ with $a \in U(\mathfrak{sl}(2,\C)$ has the obvious meaning for $\t$ extends to an isomorphism $U(\mathfrak{sl}(2,\C))\tn U(\mathfrak{sl}(2,\C)) \to U(L_\C)$ of complex algebras. As usual if $\g$ is a Lie algebra then $U(\g)$ denotes its universal enveloping algebra.

The decomposition of $L_\C$ just considered permits us to describe all infinitesimal R-matrices in $L_\C$. As usual let $t$ be the infinitesimal R-matrix in $\mathfrak{sl}(2,\C)$ coming from minus the Cartan-Killing form. Any infinitesimal R-matrix in $L_\C$ is of the form $t=a t^l+b t^r$ for some complex numbers $a$ and $b$. Some interesting combinations are following:
\begin{equation}
t_{L}=t^l-t^r=\frac{i}{8}(\s_X\tn B_X+B_X \tn \s_X+\s_Y\tn B_Y+B_Y \tn \s_Y+\s_Z\tn B_Z+B_Z \tn \s_Z)
\end{equation}
and
\begin{equation}
\hat{t}_L=t^l+t^r=\frac{1}{8} (\s_X \tn \s_X -B_X \tn B_X+\s_Y \tn \s_Y -B_Y \tn B_Y+\s_Z \tn \s_Z -B_Z \tn B_Z).
\end{equation}
They are identified with invariant non-degenerate bilinear forms in the Lie algebra of the Lorentz group. In fact, it is possible to prove that the Chern-Simons functional
\begin{equation}
CS(A)=\exp\left(i\frac{1}{4 \pi}\int_{M} \tr\left ( A \wedge dA+\frac{2}{3} A \wedge A \wedge A\right) \right )
\end{equation}
is gauge invariant if  and only if $\tr$ is defined out of the non-degenerate, invariant, bilinear form in $L$ associated with $n \hat{t}_L+s t_L$, with $n\in \Int$ and $s \in \C$, cf \cite{W}. Here $A$ denotes an $L$-valued $1$-form in a $3$ manifold $M$.

Heuristically considering the case $n\in \Int$ and $s=0$ corresponds to working with $U_q(\mathfrak{su}(2))\tn U_q(\mathfrak{su}(2))$, with $q$ chosen to be a root of unity. In our discussion, we will stick to the infinitesimal R-matrix $t_L=t^l-t^r$. The general description is not more difficult. This will correspond to working with the quantum group $U_q(\mathfrak{su}(2))\tn U_{q^{-1}}(\mathfrak{su}(2))$, or with the quantum Lorentz group $\mathcal{D}=U_q(\mathfrak{sl}(2,\C)_{\R})$ as defined by Podles and Woronowicz in \cite{PW}. We are using the notation of Buffenoir and Roche in  \cite{BR}. The knot theory obtained from $t_L$ and the unitary representations of the Lorentz Group ought to be related to any knot invariants that can be defined from the unitary representations of $U_q(\mathfrak{sl}(2,\C)_{\R}),q \in (0,1)$. These representations were originally classified by Pusz in \cite{P}. See \cite{FM} for more details.

\subsubsection{Lorentz polynomial}
Consider the infinitesimal R-matrix $t_L$ in $L_\C$ given by $t_L=t^l-t^r$. It is possible to prove that, cf \cite{GN}:
\begin{equation}
(\f_{t^l-t^r}\circ \Z) (K)=(\f_t\circ \Z)(K) \tn (\f_{-t} \circ \Z)(K).
\end{equation}
this result is a consequence of the fact $\Z(K)$ is always a group like element for any framed knot $K$, cf \cite{BN}. Also in general $ (\f_{-t} \circ \Z)(K)=(\f_{t} \circ \Z)(K^*)$. Here $K^*$ denotes the mirror image of the knot $K$.
Given $z$ and $w$ in $\C$ it is thus natural to define the Lorentz polynomial as being:
\begin{equation}
\frac{L^{z,w}}{(2z+1)(2w+1)}(K)=\frac{J^{z}}{2z+1}(K)\frac{J^{w}}{2w+1}(K^*).
\end{equation}
Here $K$ is a framed knot and $K^*$ denotes its mirror image.

Similarly with the $SL(2,\R)$ case, we have:
\begin{Theorem}
Let $R$ be a irreducible unitary representation of $SL(2,\C)$ in the complex Hilbert space $V$. There exist $z,w \in \C$ such that:
\begin{equation}
I(SL(2,\C),t_L,R)=\frac{L^{z,w}}{(2z+1)(2w+1)}.
\end{equation}
\end{Theorem}

As an example let us consider the principal series of unitary representations of $SL(2,\C)$. Some good references are \cite{GGV} and \cite {GMS}. The last reference contains an infinitesimal description of the principal series. For a more geometric description these representations in terms of hyperbolic geometry we refer to \cite{GGV} as well. The unitary principal series is parametrised by  a pair of complex numbers $z$ and $w$ with $m=z-w \in \Int$ and $i\r=z+w+1 \in i\R$. In general the parameters $m$ and $\r$ are referred to as the minimal spin and the mass of the representation. The ones of minimal spin $0$ are the balanced representations of \cite{BC}. As observed in the same reference, they admit a natural spin network theory. The representation space $V$ for $R^{m,\r}$ is $L^2(\C)$ and the action of $SL(2,\C)$ has the form:
\begin{equation}
R^{z,w}\left (\begin{pmatrix} a & b \\ c & d \end{pmatrix} \right )(f)( \x)= (b \x+d)^{z-1}(\bar{b}\bar{\x}+\bar{d})^{w-1}f\left ( \frac{a \x +c}{b \x +d}\right).\end{equation}

The representations in the principal series are unitary and irreducible, therefore there exists a knot invariant attached to them. As before, these knot invariants are particular cases of the Lorentz polynomial, in fact:
\begin{Theorem} Let $\r=z+w+1$  and $m=z-w$ we have 
\begin{equation}
I(SL(2,\C),t_L,R^{m,\r})=\frac{L^{z,w}}{(2z+1)(2w+1)}.
\end{equation}
\end{Theorem}
For a proof see \cite{FM}.

\section{Convergence Issues}

We now look at the analytic properties of the $z$-coloured Jones polynomial. As we pointed out in the introduction the power series coming out of it are in general not convergent. We now show this is what happens at least in the case of torus knots. Later we will have a look at the properties of the $z$-coloured Jones polynomial under Borel re-summation.

\subsection{On the Divergence of the $z$-coloured Jones polynomial Power Series for Torus Knots}

Let $m$ and $ p$ be two coprime positive integers. In what follows $K_{m,p}$ denotes the $(m,p)$-Torus Knot. For each $z \in \C$, Consider the following meromorphic function:
\begin{equation}
F_{m,p,z}(x)=\frac{\sinh \left ((2z+1)\sqrt{{mp}}x\right) \sinh \left (\sqrt{\frac{m}{p}}x\right)\sinh \left (\sqrt{\frac{p}{m}}x\right)}{(2z+1)\sinh \left (\sqrt{{mp}}x\right)}.
\end{equation}
Notice it is well defined if $2z+1=0$. It is an even function in $x$. Suppose $2z+1 \in \Int \setminus \{0\}$, it is possible to prove that: cf \cite{KT}
\begin{equation} \label{Kashaev}
\frac{J^{z}}{2z+1}(K_{m,p})(h)=\frac{1}{2\sqrt{\pi}}\frac{e^{-\frac{h}{4}\left (\frac{p}{m}+\frac{m}{p} \right )}}{\sinh\left( \frac{h}{2} \right )}\int_{-\infty}^{+\infty} e^{-x^2} F_{m,p,z}(\sqrt{h}x)dx, \forall h \in \C.
\end{equation}
Notice that for $2z+1 \in \Int \setminus \{0\}$ the power series $\frac{J^{z}}{2z+1}(K_{m,p})(h)$ as an infinite radius of convergence.
The origin is never a singular point of $F_{m,p,z}(x)$ for any $z \in \C$, thus if $\sqrt{h}x$ is small enough:
\begin{equation}
F_{m,p,z}(\sqrt{h}x)=\sum_{k=1} ^{+\infty} Q_{m,p,z}(k)h^kx^{2k}.
\end{equation}
In the case $2z+1 \in \Int \setminus \{0\}$ the function $F_{m,p,z}(x)$ is an entire function of exponential order. That is $f$ is analytic in $\C$ and there exist positive constants $A$ and $C$ such that $|F_{m,p,z}(x)|<Ae^{C|x|}, \forall x \in \C$. This implies that  we have the bound $ Q_{m,p,z}(k)\leq C^k/k!, k\in \{1,2,...\}$. In practice this means that, for any $h\in \C$, we can interchange the infinite summation with the integral sign in the following expression:

\begin{equation} \label{Kashaevseries}
\frac{J^{z}}{2z+1}(K_{m,p})(h)=\frac{1}{2\sqrt{\pi}}\frac{e^{-\frac{h}{4}\left (\frac{p}{m}+\frac{m}{p} \right )}}{\sinh\left( \frac{h}{2} \right )}\int_{-\infty}^{+\infty} e^{-x^2}\left ( \sum_{k=1} ^{+\infty}  Q_{m,p,z}(k)h^kx^{2k} \right )dx,\forall h\in \C,
\end{equation}
valid if  $2z+1 \in \Int \setminus \{0\}$. Therefore: 
\begin{equation}\label{Jones}
\frac{J^{z}}{2z+1}(K_{m,p})(h)=\frac{1}{2\sqrt{\pi}}\frac{e^{-\frac{h}{4}\left (\frac{p}{m}+\frac{m}{p} \right )}}{\sinh\left( \frac{h}{2} \right )}\sum_{k=1}^{+\infty}\Gamma(k+\frac{1}{2})  Q_{m,p,z}(k) h^k,
\end{equation}
if  $2z+1 \in \Int \setminus \{0\}$ and $h \in \C$. Given that the coefficients of the Taylor decomposition of
$\sinh((2z+1)\sqrt{mp}x)$ at $x=0$ depend polynomially in $z$ we show that the dependence of the coefficients $Q_{m,p,z}(k)$ in $z$ is in fact polynomial. Thus the expansion (\ref{Jones}) is true for any $z \in \C$,  now only at the power series level.   We have shown:
\begin{Theorem}Let $m$ and $p$ be coprime positive integers and $K_{m,p}$ be the $(m,p)$-torus Knots. The expansion (\ref{Jones})  of the $z$-coloured Jones polynomial of $K_{m,p}$  is correct for any $z \in \C$, as formal power series\end{Theorem}   
The functions $F_{m,p,z}(x)$ have non removable singularities if $2z+1 \notin \Int \setminus \{0\}$. In particular there exists a positive constant $C$ such that $ Q_{m,p,z}(k)> C^k$ for infinite $k$'s.  As a consequence we can conclude:
\begin{Corollary}Let $K_{m,p}$ be the $m,n$-torus knot. The power series $\frac {J^{z}}{2z+1}(K_{m,p})$ has a zero radius of convergence if $2z+1 \notin \Int \setminus \{0\}$. \end{Corollary}
It is  easy to conclude that the result above is true also for the mirror images of the class of  torus knots just considered.

It is possible to prove a similar divergence result in the Lorentz group case:
Consider the function
\begin{equation}
G_{m,p,z,w}(x)=\int_{0}^{2\pi}F_{m,p,z}(x\cos(\theta))F_{m,p,w}(ix\sin(\theta))d\theta,
\end{equation}
thus if $2z+1,2w+1\in \Int \setminus \{0\}$ we have:
\begin{equation}
\frac{L^{z,w}}{(2z+1)(2w+1)}=\frac{1}{\sinh\left( \frac{-h}{2} \right )\sinh\left( \frac{h}{2} \right )}\int_{0}^{+\infty}xG_{m,p,z,w}(x\sqrt{h})dx.
\end{equation}
If $2z+1$ and $2w+1$ are non zero integers, then $G_{m,p,z,w}$ is entire of exponential order. Otherwise it has non removable singularities. For, suppose, for example, that  $2w+1 \notin \Int \setminus \{0\}$, let $w_0$ be the first singularity of $F_{m,p,w}(ix)$ in the positive real line. It is a pole. An explicit calculation tells us  that as $x$ approaches $w_0$ from below along the real line then the first derivative of $G_{m,p,z,w}(x)$ tends to $\infty$.

Now, $G_{m,p,z,w}(x)$ is an even function with a zero of order $4$ at the origin. Put \begin{equation}G_{m,p,z,w}(x)=\sum_{k=2}^{\infty}P_{m,p,z,w}(k)x^{2k}. \end{equation}
For fixed $k$, the dependence of $P_{m,p,z,w}(k)$ in $z$ and $w$ is polynomial. In fact, if we fix $m,p,z,w$, then for $x$ small enough the series  
\begin{equation}
F_{m,p,x}(x\cos\theta)F_{m,p,w}(ix\sin(\theta))=\sum_{a,b=0}^{\infty} Q_{m,p,z}(a)Q_{m,p,w}(b)x^{a+b}{i^b}\cos^a(\theta) \sin^b(\theta)
\end{equation}
converges uniformly for $ \theta \in [0,2\pi]$. Thus
\begin{equation} 
G_{m,p,z,w}(x)=\sum_{a,b} C_{a,b}Q_{m,p,z}(a)Q_{m,p,w}(b)x^{a+b}
\end{equation}
where
\begin{equation}
 C_{a,b}=\int_{0}^{2\pi}i^b \cos^a(\theta) \sin^b(\theta)d\theta.
\end{equation}
This proves $P_{m,p,z,w}(k)$ is a polynomial in $z$ and $w$ for any $k$. Similarly as above, we conclude:
\begin{equation}
\frac{L^{z,w}}{(2z+1)(2w+1)}=\frac{1}{8\pi \sinh\left( \frac{h}{2} \right )\sinh\left( \frac{-h}{2} \right) } \sum_{k=1}^{+\infty}P_{m,p,z,w}(k) k!h^{k+1},
\end{equation}
from which follows the non-convergence of the power series defined by the Lorentz polynomial if $2z+1$ or $2w+1$ do not belong to $\Int \setminus \{0\}$,  in the case of torus Knots.

We are mainly interested in knot invariants with values in $\C$. As before we refer to them as numerical knot invariants. It is natural now to ask whether we can canonically find analytic functions of which the $z$-coloured Jones polynomial in a knot are asymptotic developments. In the remainder of this section we investigate this question in the framework of Borel re-summation.

\subsection{Borel Re-summation of Power Series}

We make now a brief description of the Borel process of re-summation of power series. We refer to \cite {SS} for full details. The paper \cite{FL} contains a simple and illuminating introduction to  this subject.

\subsubsection{Asymptotic Power Series Developments and a  Lemma Due to Borel}
Consider a sector  $\W=\{z \in \C:\theta_1 \leq \mathrm{arg}(z) \leq \theta_2, |z|\leq C\}$ in the complex plane (or in the Riemann surface of the Logarithm, if $\theta_2 -\theta_1\ge 2\pi$).  Let $\sum_{n=0}^{+\infty} a_n h^n$ be a formal power series which we do not assume to be convergent. Suppose $f:\W \to \C$ is a continuous function in $\W$ and analytic in the interior $\mathrm{int}(\W)$ of $\W$. We say that  $\sum_{n=0}^{+\infty} a_n h^n$ is an asymptotic expansion of $f$ at the origin if, for any proper subsector $\W'$ of $\W$, and any $N=1,2...$, there exists a positive constant $C_N$ such that
\begin{equation}
\left |f(h)-\sum_{n=0}^{N}a_n h^n\right | \leq C_N|h|^{N+1},\forall h \in \W'.
\end{equation}
For example the series $\sum_{n=0}^{+\infty}(-1)^n(n-1)!h^n$ is an asymptotic series for  the function $h \mapsto \int_{0}^{+\infty}e^{-\frac{x}{h}}.\frac{1}{x+1}dx$ in any sector $\W=\{h \in \C: \mathrm{Re}(h)\geq 0,|h|\leq C\}$. Here $C$ is a positive constant.

Let $\W$ be a sector in the complex plane. Denote by $G(\W)$ the space of holomorphic function in $\mathrm{int}(\W)$ and continuous in $\W$ which admit a power series asymptotic expansion at the origin. It can be proved that $G(\W)$ is closed under differentiation and it is an algebra with respect to multiplication of complex functions. If $\sum_{n=0}^{+\infty} a_n h^n$ is a asymptotic power series expansion of some $f \in G(\W)$ then the coefficients $a_n$ are unique and  can be computed by the formula:
\begin{equation}
a_n=\lim_{\substack{h \to 0\\h \in \W}}\frac{1}{n!}f^{(n)}(h).
\end{equation}
This permits us to define a map $As:G(\W) \to \C[[h]]$. It can be proved that $As$ is an algebra morphism and it is well behaved with respect to differentiation. Moreover it is a surjective map and its kernel is the space of functions analytic in $\mathrm{int}(\W)$ and continuous in $\W$ that tend  to zero faster than any $z^n$ for any $n \in \N$. This fact is known as Borel's lemma. It permits us to define a re-summation operator of formal powers series up to functions that go to zero faster that any natural power of $z$ as $z \to 0$. In practice, to re-sum a power series means finding an analytic function this power series is an asymptotic expansion of. For practical applications however it is important to reduce as much as possible the uncertainty in the process of re-summation. This is the subject of the following paragraph.

\subsubsection{Power Series in the First Gevrey Class and Formal Borel Transforms}
Let $\sum_{n=0}^{+\infty}a_n h^n$ be a formal power series. We say that it is of type Gevrey $1$ if there exists a positive constant $C$ such that $|a_n| \leq C^n n!, \forall n\in \N$. We denote by $G_1[[h]]$ the algebra of formal power series of type Gevrey $1$. If $\W$ is a sector in the complex plane, then a function in $G(\W)$ is said to be of type Gevrey $1$ if its asymptotic expansion at zero belongs to $G_1[[h]]$. It is possible to prove $G_1(\W)$ it is an algebra stable under differentiation.

Suppose the opening of the sector $\W$ is less than $\pi$. Then as  before $As_1:G_1(\W) \to G_1[[h]]$ is a surjective algebra morphism and its kernel is the space of exponential decreasing functions in $\W$. That is, functions that satisfy the estimate $f(h)\leq Ae^{-\frac{B}{|h|}}$ for some positive $A$ and $B$ in any proper subsector of $\W$. 

Recall that the Borel transform $\B(f)$ of an analytic function $f$, if it exists, is the inverse Laplace transform of it. We refer again to \cite {SS} for a discussion of this subject in the generality we need. It is well known that  $\B(1)=\d(\x)$ and $\B(h^{n+1})=\x^n/n!$. Notice we do a change of variables $h \mapsto 1/h$ in the domain of the Laplace transform. 

Let $\sum_n a_n h^n$ be a formal power series in the first Gevrey class.  Consider the generalised function 
\begin{equation}
F(\x)=\B(\sum_{n=0}^{+\infty}a_n h^n)(\x)=\sum_{n=0}^{+\infty}\frac{\a_{n+1}}{n!}\x^n+a_0 \d (\x).
\end{equation}
It is called the Formal Borel transform of $\sum_{n=0}^{+\infty}a_n h^n$. Let $0<a<C$ be a real number. Choose a direction $e^{i\theta}$ in the complex plane such that $\mathrm{Re} (e^{i\theta}/h)>0,\forall h \in \W$. Recall $\W$ opens less than $\pi$.  It is possible to prove that the incomplete Laplace Transform in the direction $e^{i\theta}$:
\begin{equation}\label{incomplete}
h \mapsto \L^{a}_{e^{i\theta}}(F(\x))(h)=\int_{0}^{ae^{i\theta}}e^{\frac{-\x}{h}}F(\x)d\x
\end{equation}
has $\sum_{n=0}^{+\infty}a_n h^n$ as a power series asymptotic development at the origin. This defines a re-summation of power series of Gevrey type 1, up to exponentially decreasing functions.

\subsubsection{Re-summation Operators}

For a better definition of the re-summation operator, we are led to consider formal power series whose Formal Borel transform can be analytically continued to the neighbourhood of some ray $e^{i\theta}\R^+_0$ in the complex plane, cf \cite{SS}. This will lead to the definition of a re-summation operator up to rapidly decreasing functions. That is, functions that verify  $\forall A > 0, \exists B >0 :|f(h)| \leq Be^{-A/|h|},\forall h \in \W'$, in any proper subsector $\W'$ of $\W$.

 Let $F(\x)$ be an analytic function possibly containing a $\d(\x)$ term, in an open subset of $\C$ containing some ray $e^{i\theta}\R^+_0$. Suppose it grows at  most exponentially along it. The Laplace transform of $F(\x)$ in the direction $e^{i\theta}$ is defined as:
\begin{equation}
\L_\ei(F(\x))(h)=\int_{0}^{+\infty e^{i\theta}}e^{-\x/h} F(\x) d\x,
\end{equation}
whenever the integral is convergent. A lot of properties of this changed Laplace transform can be deduced from the corresponding ones  of the usual Laplace Transform. In particular $\L_{\ei}(\x^{n})=\Gamma(n+1)h^{n+1},n=0,1,...$. If $a$ is a  real number bigger than $-1$ this generalises to  $\L_{1}(\x^a)=\Gamma(a+1)h^{a+1}$.

\begin{Definition} \label{resummation} ({\bf re-summation operators}) Let $\ei$ be a direction in the complex plane and $\sum_{n=0}^{+\infty} a_n h^n$ a formal power series of Gevrey type $1$. We say it is $\ei$-summable if its Formal Borel transform can be analytically continued to the neighbourhood of some ray $e^{i\theta}\R_0^+$, and if its Laplace transform converges in some nonempty subset of $\C$. If $h$ belongs to the domain of definition of the Laplace transform then the $e^{i\theta}$ re-summation operator is defined as:
\begin{equation}
S\left (e^{i\theta}, \sum_{n=0}^{+\infty}a_n h^n \right)=\L_{e^{i \theta }}\left (\B(\sum_{n=0}^{+\infty}a_n h^n)(\x)\right)(h).
\end{equation}
\end{Definition}  
We can obviously define the re-summation operators along any curve that tends to the point at infinity. 

To avoid working with generalised functions, we redefine the re-summation operators and consider:
\begin{equation}
\S\left (\ei,\sum_{n=0}^{+\infty}a_n h^n \right )=\frac{1}{h}S\left (\ei,\sum_{n=0}^{+\infty} a_n h^{n+1} \right ).
\end{equation}

\begin{Definition}{ \bf (Borel re-summability)}
A formal power series $\sum_n a_n h^n$ is said to be Borel re-summable if the re-summation operators $\S$ make sense when applied to it.
\end{Definition}

\subsection{Back to Knots!}

Let $K$ be any framed knot. Consider a unitary representation $R$ of $SL(2,\R)$. It is possible to obtain some estimates for the coefficients of the Kontsevich Integral as well as for the matrix elements of $R^{\infty}$ that are suitable to prove that $I(SL(2,\R),t,R)$ defines a formal power series of Gevrey type $1$. In fact: 

\begin{Theorem} \label{Gevrey} Let $K$ be a framed Knot, $z$ and $w$ two complex numbers. The series defined by $\frac{J^z}{2z+1}(K)$ and $\frac{L^{z,w}}{(2z+1)(2w+1)}(K)$ are of Gevrey type $1$. That is, their Formal Borel Transform defines an analytic function in a neighbourhood of zero.  
\end{Theorem}

We will sketch a proof of this theorem in the appendix. It is a purely technical proof and the rest of the article is fairly independent of this result.

\subsubsection{The Case of Torus Knots}
We now analyse in detail the properties of the Formal Borel Transform of $\Jz(K)$ in the case $K$ is a torus knot. We shall see the re-summation operators make sense in this case. The ambiguity in the process of re-summation is reduced to a finite set of functions differing by rapidly decreasing functions.

Let  $f$ and $g$ be two complex valued functions continuous  in some open subset of $\C$ which contains $0$. Recall that their convolution is defined as
\begin{equation}
(f*g)(x)=\int_{0}^{x} f(x-\x)g(\x) d\x, \end{equation}
whenever it makes sense. The contour of integration is chosen to be the segment connecting $0$ and $x$. If $f$ and $g$ are entire then their convolution is an entire function. It is well know that if $f$ and $g$ grow at most exponentially along a ray $e^{i\theta}$ the the same is true for their convolution, and in addition $\L_{i\theta}(f*g)(h)=\L_{i\theta}(f)(h)\L_{i\theta}(g)(h)$ whenever all transforms makes sense. Suppose $f$ and $g$ are meromorphic functions which do not have the  origin as a singular point. Let $D(f,g)$ be the set of points in $\C$ that can be connected to the origin by a straight line which does not pass by any singularity of $f$ or $g$. The convolution of $f$ and $g$ is an analytic function in $D(f,g)$. The singularities of $f*g$ are in general ramifying due to the fact $f$ and $g$ may have non zero residua at their singular points. One can easily determine the Taylor series of $f*g$ in zero from the Taylor series of $f$ and $g$, in fact if $f(x)=\sum_{n=0}^{+\infty} f_n x^n$ and $g(x)=\sum_{n=0}^{\infty}$ then $(f*g)(x)=\sum_{n,m =0}^{\infty}\frac{m!n!}{(m+n+1)!}x^{m+n+1}$, which converges whenever the Taylor series of $f$ and $g$ converge. This is because,  in general, $x^m * x^n =\frac{m!n!}{(m+n+1)!}x^{m+n+1};m,n \in \N$,  apart from some simple analysis,. In particular if the coefficients of the Taylor series of $f_z(x)$ and $g_z(x)$ in zero depend polynomially in $z$ then the same is true for the coefficients of the Taylor series of $f*g$.

 Suppose $f$ is an odd meromorphic function non singular at the origin. Then the convolution $1/\sqrt{x} * f(\sqrt{x})$ defines an analytic function in $D(f)$, where $D(f)$ is the set of points in the complex plane that  can be connected to zero by a straight line whose square does not meet any singularity of $f$. Notice $D(f)$ contains a neighbourhood of $0$. In fact if $f(x)=\sum_{n=0}^{+\infty}a_n x^{2n+1}$ in a neighbourhood of the origin then $1/\sqrt{x} * f(\sqrt{x})=\sum_{n=1}^\infty \frac{\Gamma(1/2)\Gamma(n+1+1/2)}{\Gamma(n+2)}a_n x^n$, in the same neighbourhood. As before, the singular points of $1/\sqrt(x) * f(\sqrt{x})$ are usually ramifying. Consider
\begin{align}
I_{m,p,z}
        &=\frac{\Gamma(\frac{1}{2})^{-1}}{\sqrt{x}}*\frac{F_{m,p,z}(\sqrt{x})}{2\sqrt{x}}.
\end{align}
Thus if $2z+1\in \Int \setminus \{0\}$ then $I_{m,p,z}(x)$ is an integral
function of exponential order. In general for any $z \in \C$, the function
$I_{m,p,z}(x)$ is meromorphic, but it is well defined in  a neighbourhood of
$0$. In addition the Taylor coefficients of $I_{m,p,z}$ depend polynomially in
$Z$, since the same is true for the Taylor coefficients of  $F_{m,p,z}(x)$ in
$0$, fact we have seen before. Let
\begin{equation}
 H_{m,p,z}(x)=\B\left (\frac{h}{2\sqrt{\pi}}\frac{e^{-\frac{h}{4}\left (\frac{p}{m}+\frac{m}{p} \right )}}{\sinh\left( \frac{h}{2} \right )} \right ).
\end{equation}
It is an entire function of $x$. In fact it is a function of exponential order, the order of which is smaller than $1/\pi$. The Laplace transform in any direction  extends analytically back the original function. Suppose $m$ and $p$ be coprime positive integers and let $K_{m,p}$ denote the $(m,p)$-Torus Knot. From equation (\ref{Kashaev}) it is imediate that if $2z+1\in \Int \setminus \{0\} :$ we have 
\begin{equation}\label{BorelofJones}
\B(h \frac{J^z}{2z+1}(K_{m,p}))=H_{m,p,z}(x) *I_{m,p,z}(x), \forall x \in \C.
\end{equation}
The proof of this fact uses only elementary theory of Laplace transforms. Now $H_{m,p,z}(x) *I_{m,p,z}(x)$ is a $z$-dependent family of analytic functions in a neighbourhood of zero,  whose Taylor series coefficients in zero depend polynomially in $z$. In particular, equation (\ref{BorelofJones}) is correct for any $z \in C$, for $x$ in a neighbourhood of $0$, thus always in the power series level.
Suppose  $2z+1$ is not a non zero integer. Then  $\B(h \frac{J^z}{2z+1}(K_{m,p}))(x)$ is an analytic function in the set $\C\setminus (-\infty,\pi^2/(mp)]$. It is possible  to prove that the  singularity in $-\pi^2/(mp)$ is ramifying, thus we can not remove it. This tells us that we cannot refine the estimate in theorem \ref{Gevrey}.

We summarise our discussion and add some more simple facts:
\begin{Theorem}\label{Borel}
Let $m$ and $p$ be coprime positive integers and $K_{m,p}$ be the $(m,p)$-Torus Knot. Let $z$ be a complex number such that $2z+1 \notin \Int \setminus\{0\}$ The Formal Borel transform of $h \frac{J^z}{2z+1}(K_{m,p})$ extends to an analytic  function in  $\C\setminus(-\infty,\pi^2/(mp)]$. The analytic extension follows from (\ref{BorelofJones}). The singularity in  $-\pi^2/(mp)$ is a ramification point. If $\ei\R^+_0$ is a ray in the complex plane that is not the set of non-positive real numbers then  $\B(h \frac{J^z}{2z+1}(K_{m,p}))(x)$ grows not faster than $Ae^{|x|/\pi}$ along it, and the constant $A$ is independent of the direction chosen.
\end{Theorem}

Recall now definition \ref{resummation} and the comments after. Let $B(\pi)=\{h \in \C, |h| < \pi$. If $h \in B(\pi)$ consider $D(h)=\{w \in {S^1 \setminus \{-1\}}: \mathrm{Re}(w/h)>1/\pi\} \neq \emptyset$. Given $h\in B(\pi)$ then $D(h)$ is connected if $\mathrm{Re} (h) \geq 0$ and it has two connected components otherwise. 

Suppose $\mathrm{Re}(h)\geq 0$.  Given that $D(h)$ only has one connected component, the re-summation  $\S(w,h\Jz(K_{m,p})(h))$ does not depend on $w \in D(h)$. Also, if $2z+1$ is a real  or an imaginary number, it takes positive real values of $h$ to real numbers. Notice the last case corresponds to the representations of $SL(2,\R)$ in the principal series. Due to the fact that $D(x)$ has two connected components if $\mathrm{Re}(x)<0$ the same is not true in this case. In fact if $x$ is a real number smaller than $0$ the two different re-summations will give conjugated complex numbers of non zero imaginary part, for the case $2z+1$ is real or imaginary. Notice that the values of the re-summation, as well as the possible ambiguities can be read directly from equation (\ref{Kashaev}). 

The case of the mirror images of torus knots considered can be treated in a similar way. In this case the domain of re-summation is also $B(\pi)$. The re-summation procedure having $2$ branches if $\mathrm{Re}(h)>0$ and one otherwise. Like before, the real  parts of the re-summation for real $h$, are independent of the re-summation procedure if we stick to representations in the principal series. In the Lorentz group case, the domain for the re-summation of $\frac {L^{w,z}}{(2z+1)(2w+1)}(K_{m,p})$ is again  $B(\pi)$. Unless $h$ is imaginary,  we now have always two branches for re-summation, if $2z+1$ and $2w+1$ are not non zero integers. Therefore
\begin{Theorem} Let $K_{m,p}$ be the $(m,p)$-torus knot. For any $z,w \in \C$ the power series $\Jz(K_{m,p})(h)$ and  $\frac {L^{w,z}}{(2z+1)(2w+1)}(K_{m,p})(h)$ are Borel re-summable.
\end{Theorem}

\subsection{Conclusion} 
Let $G$ be any real Lie group. As we have seen, if we are given an $\ad$-invariant non degenerate symmetric bilinear form in $\g \tn_\R \C$, there exists a framed knot invariant $\f_t \circ \Z$ with values in the algebra of formal powers series in $h$ over the centre of the universal enveloping algebra of $\g \tn_\R \C$, in other words is is  a formal power series of differential operators in $G$. If $\g \tn_\R \C$ is a semisimple Lie algebra and $V$ is a finite dimensional representation of $G$ then this formal power series at a particular value of $h$ evaluates to an operator, or $h$-dependent family of operators, $V \to V$, which is a $G$-intertwiner. The trace of these operators will yield the usual quantum group knot invariants.  

Instead of trying to make sense of this power series of differential operators at a particular value of $h$ for the case $V$ is an infinite dimensional representation, it is observed that each term of this power series acts as a multiple of the identity operator in the space of smooth vectors of any unitary irreducible representation. Therefore $\f_t \circ \Z$ can be evaluated in any representation of this type. This gives us a knot invariant with values in the algebra of formal power series over $\C$ for any unitary irreducible representation of $G$. At least in the case $\g$ is semisimple these knot invariants express out of the knot invariants associated with finite dimensional representations of the quantised universal enveloping algebra $U_h(\g)$ by a  analytic continuation in the power series level, see theorem \ref{Continuation} and the comments after. Notice that a coherent way of realising this analytic continuation at the level of the values of the quantum group invariants would give us a way of defining $\C$-valued knot invariants associated with infinite dimensional unitary representations. 

However, at least in the case $G=SL(2,\R)$ or $G=SL(2,\C)$, for some interesting infinite dimensional unitary representations the power series that we obtain have zero radius of convergence. The process of Borel re-summation is analysed in the case of the $z$-coloured Jones polynomial. It is observed that in this case the power series are of Gevrey type $1$, which means their Formal Borel Transform defines an analytic function in a neighbourhood of zero. In particular we can define a re-summation operator apart from exponentially decreasing terms, for example from expression (\ref{incomplete}). However for the Borel re-summation procedure to be applied in full generality, we need the extra conditions that the formal Borel transform can be analytically continued to some ray $\ei \R_0^+$  admitting a Laplace transform with a non zero domain along this direction. This reduces the uncertainty in the re-summation to a numerable set of functions differing by rapidly decreasing terms. It is verified that this is indeed what happens when considering torus knots. Namely the formal Borel transform can be continued along mostly all directions. The singularities of it are ramifying, and the result of re-summation can have at most two branches, both  solutions of the problem of analytically extending the values of the coloured Jones polynomial to complex spins. In general for representations in the principal series with real $h$ the two branches define conjugated complex numbers. The question that arises naturally concerns the behaviour of the $z$-coloured Jones polynomial in the general case, namely:
  
\begin{enumerate}
\item Can the Formal Borel Transform of $\frac{J^{z}}{2z+1}(K)$ be analytically continued in general?

\item Does it admit Laplace transforms with non zero domain?
\item Which kind of ambiguities will result from this process of re-summation?
\item Do these ambiguities have a meaning?
\end{enumerate}

\section*{Acknowledgements}
This work was realised in the course of my PhD program in the University of Nottingham under the supervision of Dr John W. Barrett. I was financially supported by the programme {\em `` PRAXIS-XXI''}, grant number SFRH/BD/1004/2000
of {\em Funda\c{c}\~{a}o para a Ci\^{e}ncia e a Tecnologia} (FCT), financed by the  European Community fund {\em Quadro Comunitario de Apoio III}, and also by {\em Programa Operacional ``Ci\^{e}ncia, Tecnologia, Inova\c{c}\~{a}o''} (POCTI) of the {\em Funda\c{c}\~{a}o para a Ci\^{e}ncia e a Tecnologia } (FCT), cofinanced by the European Community fund FEDER.

\section*{Apendix (Proof of Theorem \ref{Gevrey})} 

We assume that the reader is familiar with the construction of the Kontsevich
Integral, as well as the algebraic structure in the space of chord diagrams.
See for example \cite{BN}, \cite{SW} or \cite{CS}. Let $\A_\fin=\oplus_{n \in \N_0}
\A_n$. The connected sum of chord diagrams defines a product  $\A_n \tn \A_m
\to \A_{m+n}$, providing $\A_\fin$ a graded algebra structure. By definition,
the algebra of chord diagrams $\A$ is the graded completion of $\A_\fin$, or
alternatively the algebra of formal linear combinations of the form $\sum_n
w_n h^n$ where $w_n \in \A_n, \forall n \in \N_0$.

\begin{figure}
\includegraphics{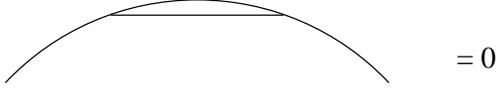}
\caption{\label{framingindependence} Framing independence relation on chord diagrams}
\end{figure}

 The unframed Kontsevich Integral has values on the algebra $\A'$, defined
 similarly to the algebra $\A$ of chord diagrams, but considering also the
 framing independence condition of figure
 \ref{framingindependence}. Alternatively we can define  $\A'$ as the graded
 completion of $\A'_\fin=\A_\fin/<\ominus \A\_\fin>$, the last with the
 obvious grading, see \cite{W2}. Here $\ominus$ denotes the chord diagram with
 just one chord.

\begin{figure}
\includegraphics{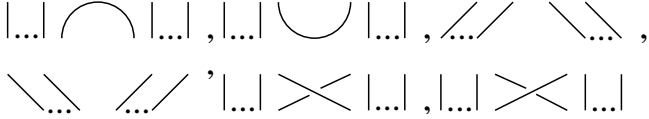}
\caption{\label{Generators} Generator tangles of the kinds $\cap$, $\cup$,
  $C_-$, $C_+$, $X_-$ and $X_+$}
\end{figure}

Consider a knot $K$ made out of $n$  tangle generators $G$ like the ones shown
  in figure \ref{Generators} on top of each others, with a chosen
  orientation. They can be of six different kinds, namely  $\cap$, $\cup$,
  $C_-$, $C_+$, $X_-$ and $X_+$. We can suppose they have height $1$. The knot $K$ is a Morse knot, and the critical points of it are contained in $\{0,...,n\}$. Let $I_k=[k-1,k]$, where $k=1,...,n$.

\begin{figure}
\includegraphics{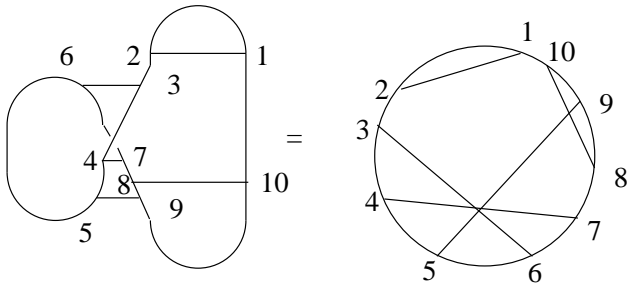}
\caption{\label{PairingChord} Defining a chord diagram out of a pairing.}
\end{figure}

Define:
\begin{multline}
Z(K)=\sum_{m =0}^{\infty}\frac{h^m}{(2\pi i)^m}
\int_{\substack{0<t_1<t_2<...<t_m<n\\t_j \notin\{0,...,n\}, j=1,...,m }}\\
\sum_{\substack{\textrm{pairings}\quad P\\P=\{\{z_j,z_j'\}:I_{k_j}\to
    \C\}_{j=1}^m}} (-1)^{\#
  \downarrow P} w_P \bigwedge_{j=1}^m\frac{dz_j-dz_j'}{z_j-z_j'}\in \A'
,
\end{multline}
that is
\begin{equation}\label{todecript}
Z(K)=\sum_{m=0}^{\infty}h^m\sum_{\substack{\textrm{pairings}\quad P\\P=\{\{z_j,z_j'\}:I_{k_j}\to
    \C\}_{j=1}^m}}Z(P,K)w_P\in \A'.
\end{equation}
Here $I_{k_j}$ denotes the interval where the pair of functions
$p_j=\{z_j,z'_j\}:I_{k_j} \to \C$  (a chord) is defined, thus
$k_1<k_2<...<k_m$. If $P$ is a pairing then $w_P$ denotes the chord diagram
constructed out of it as in figure \ref{PairingChord}, whereas $\downarrow P$
denotes the set of $z_i$ or $z'_i$ where the orientation of $K$ points
downwards. The framing independence relation ensures the integrals are
convergent, for the divergent integrals evaluate to zero in $\A'$, see
\cite{BN} $4.3.1$.

\begin{figure}
\includegraphics{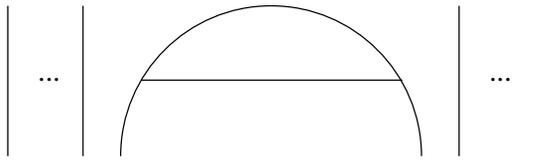}
\caption{\label{alfa} A chord of type $A$}
\end{figure}

\begin{figure}
\includegraphics{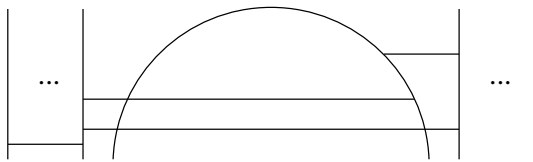}
\caption{\label{beta} Chords of type $B$}
\end{figure}

Any pairing $P$ with $m$ chords is defined in some connected component of
$\{0<t_1<t_2<...<t_m<n,t_j \notin\{0,...,n\}, j=1,...,m\}$, in particular
there are $m_k$ chords $p_i=\{z_i,z_i'\}$ defined in each $I_k$ for $k \in
\{1,...,n\}$. We can thus index the chords as
$p^k_i,k=1,....n,i=1,...m_k$. Let $\mm=(m_1,...,m_k)$, thus
$|\mm|=m_1+...+m_n=m$. We say a chord $p$ defined in a $I_k$ associated with a
generator $G_k$ of the kind $\cup$ or $\cap$ is of type $T(p)=A$ (resp
$T(p)=B)$ if it looks like the one in figure \ref{alfa} (resp the ones in
figure \ref{beta}). If the generator $G_k$ is of type $C_-$, $ C_+$, $X_-$ or $X_+$ then by definition all the chords in $I_k$ are of type $B$. If $G_k$ is a generator of type $\cup$, $X_+$ or $X_-$, let $T^k=\left (T(p^k_1),...,T(p^k_{m_k})\right )$. If $G_k$ is of type $\cap$ then  we define $T^k=\left (T(p^k_{m_k}),...,T(p^k_{1})\right )$. Let $B(k)$ be the number of chords in $I_k$ of type $B$ and $B(P)$ be the total number of $B$-chords in $P$. Notice that since we are considering the framing independence relation we can suppose $T^k_1=B,k=1,...,n$. This is a necessary condition for all the integrals to be convergent, in the first place.

Explicit calculations prove:
\begin{Lemma}\label{Intestimate} Fix a knot $K$ as before, there exists a $C<+ \infty$ such that for any pairing $P$ with $m$ chords we have:   
\begin{equation}
|Z(P,K)| \leq  C^m\prod_{k=1}^n\frac{1}{B(k)!} \prod_{k=1}^n \prod_{\substack{i\in \{1,...,m_k\}\\T^k_i=A}}\frac{1}{\# \{j\in \{1,..,i\}:T^k_j=B\}}.
\end{equation}
\end{Lemma}
\begin{Proof}{\bf (sketch)}  Let
$T=(T_1,..,T_n)$ be a sequence of $A$'s and $B$'s with $T_1=B$. Let  $B(T)$ be
the number of elements of $T$ equal to  $B$, and let $C(A)=1$ and $C(B)=1/2$. We have:
\begin{multline}
I(T)=\int_{0<t_1<...<t_m<1} \prod_{i=1}^{m}\frac{1} {t_i^{C(T(A_i))}}dt_1...dt_m=2^{m}\frac{1}{B(T)!}\\ \prod_{\substack {i\in \{1,...,m\}\\T_i=A}} \frac{1}{\# \{j \in \{1,...,i\}:T_j=B\} }.
\end{multline} 
One proves this  from the equality 
\begin{equation}
\int_{-\infty<x_1<x_2<...<x_n<0}e^{\sum_{i=1}^n\l_i x_i}dx_1...dx_n\\=\frac{1}{\l_1(\l_1+\l_2)...(\l_1+\l_2+...+\l_n)},
\end{equation}
(easy to prove by induction  if $\l_1,\l_2,...,\l_n>0$) by a continuity argument.  Consider the generators $\cup$ and $\cap$ to be made of
semicircles of radius $1$ and strings parallel to $z$-axis. Unpacking equation \ref{todecript} yields a product of $n$ integrals,
 one for
each $I_k$, equal to, or bounded by, integrals like $C^{m_k} I(T)$, where $C$
is fixed. 
\end{Proof}

Let $z\in \C$ be a complex number, and let $\l_z$ be the central character of
the representation of $\sl$ of highest weight $2z$. We wish to give an
estimate for $|\l_z(\f_s(w_P))|$, where $s$ denotes the infinitesimal
$R$-matrix in $\sl$ coming from a quarter of the Cartan-Killing form and $w_P$
is the chord diagrams made out of the pairing $P$. Even though $\f_s(w)$ does
not make sense in $\A'$, since  $\f_s$ does not satisfy the framing
independence relations, $w_P$ does define a chord diagram, thus an element of
$\A$, and in particular  it makes sense to consider $\f_s(w_P)$.  

Recall the recursive evaluation of the $\sl$ weight system $\f_s:\A \to
U(\sl)$ in \cite{CV}. Let $w$ be a chord diagram and let $a$ be a chord of it,
thus $a$ divides the circle supporting $w$ into two semicircles. Define $C(a,w)$ as the set
of chords of $w$ that cross $a$. It has cardinality $x(a,w)$.  We define
$w(a)$ as the chord diagram obtained from $w$ by removing $a$. Let  $b$ and
$c$ be distinct chords in $C(a,w)$, they have endpoints $e_a,e_b,f_a,f_b$ such
that $e_a$ and $e_b$ (resp $f_a$ and $f_b$) lie in the same semicircle. Define
$w^{X}(a,b,c)$ (resp $w^{||}(a,b,c)$) as the chord diagrams obtained from $w$ through removing $a$,$b$ and $c$ and adding two new chords joining $e_a$ with $f_b$ and $e_b$ with $f_a$ (resp $e_a$ with $e_b$ and $f_a$ with $f_b$). Theorem $1$ of \cite{CV} tells us:
\begin{multline}\label{recursive}  
(\l_z \circ \f_s)(w)=(\l_z \circ \f_s)(\ominus w(a) - 2x(a,w) (\l_z \circ \f_s)(w(a))+\\\sum_{\substack{b,c \in C(a,w)\\b \ne c} }2 \left ( (\l_z \circ \f_s)(w^{X}(a,b,c)) -(\l_z \circ \f_s)(w^{||}(a,b,c)) \right ).
\end{multline}     
Let $w$ be a chord diagram with $m$ chords having a set $\a(w)$ of chords that
do not cross each others (we do not suppose that this set is maximal). Pick up
a chord $a_0\in \a(w)$. Consider a diagram $w'$ appearing in the recursive
evaluation \ref{recursive}. Then $\a(w')=\a(w)\setminus \{a_0\}$ is a set of
chords of $w'$ that do not intersect each others. In addition we always have
$x(a,w')\le x(a,w), \forall a \in \a(w')$. We need to use the fact that the
chords of $\a(w)$ do not cross each others to prove this. An obvious induction
based on \ref{recursive} and this last fact tells you:
\begin{align}
|(\l_z \circ \f_s)|(w)&\leq S_z(m-\#\a(w))\prod_{a \in a(w)}\left( |c_z|+2x(a,w)^2\right) \\
&\leq S_z(m-\#\a(w))C_z^{\#\a(w)}\prod_{a \in a(w)}\left(1+x(a,w)\right)^2.
\end{align}
Here $c_z=(\l_z \circ \f_s)(\ominus)$ and $C_z=\max(|c_z|,2)$. In general, if $i  \in \N$, we define $S_z(i)=\max\{|\l_z \circ \f_s(x)|,x\in W_{i}\}$ where $W_{i}$ is the set of chord diagrams with $k$ chords. The recursive evaluation tells you again that $S_z(i) \le \prod_{j=1}^i (|c_z|+2 j^2)\leq C_z^j(j+1)!^2$. Given this last estimates we prove:
\begin{Lemma}\label{Mestimate}
Let $z \in \C$ there exist a constant $C_z <+\infty$ such that the following is true:  Let  $w$ be a chord diagrams with  $m$ chords  having a set  $\a(w)$ of chords that do not cross each others. We have the estimate:
\begin{equation}\label{priorMestimate}
|(\l_z \circ \f_s)|(w)<C_z^m (k+1)!^2\prod_{a \in a(w)}\left( 1+x(a,w)\right)^ 2,
\end{equation}
where $k=m-\#\a(w)$. 
If in addition $x(a,w)\ne 0, \forall a\in \a(w)$ we can write this in a more useful form for later, namely 
\begin{equation}
|(\l_z \circ \f_s)|(w)<(2C_z^m) (k+1)!^2\prod_{a \in a(w)}\left(x(a,w)\right)^ 2.
\end{equation}
\end{Lemma}
Let $K$ be an oriented knot and let $P$ be a pairing with $m$ chords. We want
to apply the lemma above to $w_P$, the chord diagram made out of $P$. Each
chord $p$ of $P$ gives rise to a chord $a_p$ of $w_P$. Two chords $a_p,a_{p'}$ with $p$ and $p'$  chords of type $A$ cannot cross each others, therefore we define $\a(w_P)=\{a_p: T(p)=A\}$. If $p$ is defined in a generator $G$ of $K$ of type $\cup$ (resp. $\cap$), then $x(a_p,w_P)$ is the number of $B$-chords defined in the same generator $G$ staying bellow (resp. above) $p$, thus $x(a_p,w_P)$ can be supposed to be different of zero by the framing independence relation, for $T^k_1=B$. Combining  lemmas \ref{Intestimate} and \ref{Mestimate}, we conclude there exists a $C<+\infty$ such that: 
\begin{equation}
|(\l_z \circ \f_s)Z(P,K)w_P| \leq C^m \frac{(B(P)!)^2}{\prod_{k=1}^nB(k)!} \prod_{k=1}^n \prod_{\substack{i\in \{1,...,m_k\}\\T^k_i=A}}{\# \{j\in \{1,..,i\}:T^k_j=B\}}
\end{equation}
where $C$ only depends on $z$ and $K$. Notice $B(P)=B(1)+...+B(n)$ is the number of chords in  $P$ of type $B$.

An easy consequence of Stirling inequalities is the fact that given $n \in \N$, there exists a $C< \infty$ such that $  m!\le C^m m_1!...m_n!$ where $m_1+...+m_n=m$. Another consequence is the fact that there exists a $C<+\infty$ such that $m^m<C^mm!,\forall m \in \N$. Let $K$ be a knot and $P$ a pairing with $m$ chords. Putting everything together we prove 
\begin{align}
|(\l_z \circ \f_s)(Z(P,K)w_P)|&\leq C^m B(P)^{B(P)} B(P)^{A(P)}\\&\leq C^m(A(P)+B(P))^{A(P)+B(P)}\\&\leq D^m m!,
\end{align}
for any pairing with $m$ chords. Here $A(P)$ is the number of chords of $P$ of type $A$.  Notice $D<+\infty$ only depends on $z$ and $K$. We know $s=-t/4$ thus $f_t(w)=(-4)^m \f_s(w)$, where $m$ is the number of chords of $w$. therefore we can also find a constant $C< +\infty$ such that:
\begin{equation}\label{Prior}
|(\l_z \circ \f_t)(Z(P,K)w_P)|\le C^mm!.
\end{equation}
We have now done the most difficult part of the proof of theorem \ref{Gevrey}, even though this last inequality is not quite enough yet. Let us go back to the definition of the Kontsevich Integral. Let $\infty$ be a Morse knot equivalent to the unknot but with $4$ critical points. Let $N$ be the number of critical points of $K$. The unframed Kontsevich integral of $K$ is
\begin{equation}
\cZ_u(K)=\frac{Z(K)}{Z(\infty)^{N/2}}.
\end{equation}
Let us define the framed version of the Kontsevich Integral. Recall we have a coproduct in $\A_\fin$ defined by $D(w)=\sum_{x \subset w} x \tn (w\setminus x)$, whenever $w$ is a chord diagram. Here $x$ is a chord diagram made out some chords of $w$ and $w \setminus x$ is the complementary diagram. Recall $\A_\fin$ has a grading $\deg$ where the grading coefficient of a chord diagram is its number of chords. Consider the map $\psi:\A_\fin \to \A_\fin$ such that 
\begin{equation}
\psi(w)=\sum_{x \subset w} (-\ominus)^{\deg(x)} (w\setminus x).
\end{equation}
for a chord diagram $w$. Thus $\psi$ satisfies the $4T$-relation. In fact $\psi$ is a Hopf algebra projection, see \cite{W2} (this is not a trivial fact). The morphism $\psi$ is zero on the ideal generated by $\ominus$, thus $\psi$ defines an algebra morphism $\A'_\fin \to \A_\fin$. It extends to a morphism $\psi_0:\A' \to \A$ of the graded completions.
\begin{Definition}
The Kontsevich Integral of $K$ is
\begin{equation}
\cZ(K)=e^{\ominus F(K)h}\psi_0(\Z_u(K))\in \A,
\end{equation}
where $F(K)$ is the framing coefficient of $K$.
\end{Definition}
This definition is equivalent to the standard one. See \cite{SW} or \cite{LM2} theorem $5.13$. Notice $\psi(Z(\infty))^{-1}=\cZ(O)$, which is the (framed) Kontsevich Integral of the zero framed unknot. Let $z\in \C$. Since $(\l_z  \circ \f_t)$ is an algebra morphism, we have:
\begin{align}
  \frac{J^z(K)}{2z+1}&=(\l_z  \circ \f_t)\cZ(K)\\ &=(\l_z  \circ \f_s)\left (e^{\ominus F(K)h}\right )  (\l_z  \circ \f_s)\left (\psi_0(Z(\infty))^{-N/2} \right ) (\l_z  \circ \f_t)(\psi_0(Z(K)))\\  &=e^{c_z h}(\l_z  \circ \f_t)\left(\cZ_f(O)\right)^{N/2}(\l_z  \circ \f_t)(\psi_0(Z(K))\\ &=e^{c_z h}\left (\frac{1}{2z+1}\frac{\sinh((2z+1)h/2)}{\sinh(h/2)}\right )^{N/2} (\l_z  \circ \f_t)(\psi_0(Z(K)).
\end{align}
Since the set $G_1\h$ of power series of Gevrey type $1$ forms an subalgebra of $\C\h$, the proof of theorem \ref{Gevrey} will be finished if we prove that
\begin{equation}\label{toprove}
(\l_z  \circ \f_t)(\psi_0(Z(K))=\sum_{m=0}^{\infty}h^m\sum_{\substack{\textrm{pairings}\quad P\\P=\{\{z_j,z_j':I_{k_j}\to \C\},j=1,...,m\}\}}}(\l_z  \circ \f_t)(\psi_0 (Z(P,K)w_P)
\end{equation}
is of Gevrey type $1$. The estimate for $|\l_z(\f_s)(w)|$ in \ref{priorMestimate} continues to hold if we remove a chord from $w$, maintaining the right hand side of \ref{priorMestimate} fixed.  Immediately we have:
\begin{align}
|(\l_z \circ \f_s)|(\psi_0(w))&\le \sum_{x \subset w}|(\l_z \circ \f_s) (-\ominus)|^{\deg(x)}| (\l_z \circ \f_s) (w\setminus x)|\\
&\le (2C_z)^m C_z^m (k+1)!^2\prod_{a \in a(w)}\left( 1+x(a,w)\right)^ 2.
\end{align}
Recall that  $c_z=(\l_z \circ \f_s)(\ominus)$ and $C_z=\max(|c_z|,1)$. Notice there are $2^m$ splittings $w=x \cup (w\setminus x)$ if $w$ is a chord diagram with $m$ chords. As before we prove.
\begin{Lemma}\label{Prior2} The estimate \ref{Prior} continues to hold if we put $(\psi_0 (Z(P,K)w_P)$ instead of $(Z(P,K)w_P)$.
\end{Lemma}
We now prove $(\l_z\circ \f_t)(\psi_0(Z(K))\in G_1\h$. Since we have proved lemma \ref{Prior2}, a glance at equation \ref{toprove} reduces the proof of this to estimating the number of pairings $P$ with $m$-chords, which is a simple exercise of counting. Let $M=\max_{t \in \R} \#((t\times \C )\cap K)$. Consider pairings  $P$ with $m$ chords, having  $m_i$ chords in each interval $I_i$ for $i=1,...,n$, thus $m_1+...m_n=m$. There are at most $(M(M-1)/2)^m$ pairings like this. Recall  the classical combinatorics problem which asks about the number 
of ways we can put $X$ indistinguishable objects into $N$ distinguishable 
boxes. Its solution is  $\frac{(N+X+1)!}{(N-1)!X!}$. In our case we have exactly $X=m$ objects (chords) and $N=n$ boxes (intervals $I_k$). Thus there are at most $(M(M-1)/2)^m\frac{(n+m+1)!}{(n-1)!m!}$ pairings $P$ with $m$ chords. Given that $n$ and $M$ are constant, this last term grows exponentially with respect to $m$. This finishes the proof of theorem \ref{Gevrey}.

\end{document}